\documentclass{amsart} 
\usepackage[mathscr]{eucal}
\usepackage{amsmath,amsfonts}
\newtheorem{theorem}{Theorem}[section]

\newtheorem{corollary}[theorem]{Corollary}

\makeatletter
\@addtoreset{equation}{section}
\makeatother

\newcommand{\cA}{{\mathcal A}}

\newcommand{\cD}{{\mathcal D}}

\newcommand{\cH}{{\mathcal H}}

\newcommand{\cL}{{\mathcal L}}
\newcommand{\cM}{{\mathcal M}}

\newcommand{\cP}{{\mathcal P}}

\newcommand{\cW}{{\mathcal W}}

\newcommand{\CC}{{\mathbb C}}
\newcommand{\NN}{{\mathbb N}}

\newcommand{\FF}{{\mathbb F}}

\begin{document}
\pagestyle{plain}

%\begin{flushright}
%  \it Date of this draft: \today
%\end{flushright}
\bigskip

\title{Orthogonal polynomials  \\
in several non-commuting variables. II} 
\author{ T. Banks }\author{T. Constantinescu} 

\address{Department of Mathematics \\
  University of Texas at Dallas \\
  Box 830688, Richardson, TX 75083-0688, U. S. A.}
\email{\tt banks@utdallas.edu} 
\email{\tt tiberiu@utdallas.edu}

\begin{abstract}
In this paper we continue to investigate a certain class of Hankel-like 
positive definite kernels using their associated orthogonal polynomials.
The main result of this paper is about the structure of this kind of kernels.

\end{abstract}

\maketitle

\section{Introduction}
Positive definite kernels are studied for their
manifold applications. In this paper we consider 
a special type of kernels $K$ defined on the free 
semigroup on $N$ generators 
with the property that 
$$K(\alpha \sigma ,\tau )=K(\sigma , I(\alpha )\tau )$$
for any words $\alpha ,\sigma ,\tau $, where $I(\alpha )$
denotes the word obtained by writting $\alpha $ in the reverse order.
These kernels appear in many situations, see for instance
\cite{CG1} and \cite{CG2}.
Our goal is to determine an explicit structure of the positive definite 
kernels satisfying the above invariance property. Since in the 
case $N=1$ such kind of kernels are precisely the Hankel kernels, 
it is quite natural to consider associated orthogonal polynomials 
and to study their properties. Our main result establishes the 
connection between moments and Jacobi coefficients, as a multivariable
extension of a classical result. We also describe the Jacobi 
coefficients of the free products of orthogonal polynomials.

\section{Orthogonal polynomials}
We introduce orthogonal polynomials on several hermitian variables and we
discuss several general results. Especially, we emphasize the usefulness 
of a matrix notation that reduces very much the degree of complexity and 
makes clear the analogy with the classical, one-dimensional case.
Let $\FF _N^+$ be the unital free semigroup 
on $N$ generators $1,\ldots ,N$ with lexicographic
order $\prec $. In particular, $\FF ^+_1$ is the set $\NN _0$ of nonnegative
integers. The set of positive integers
will be denoted by $\NN $. The empty word is the identity element
of $\FF ^+_N$ and the length of the word $\sigma $ is 
denoted by $|\sigma |$. The length of the empty word is $0$.
There is a natural involution on
$\FF _N^+$
given by $I(i_1\ldots i_l)=i_l\ldots i_1$ as well 
as a natural
action of $\FF _N^+$ on itself by juxtaposition, 
$(\sigma,\tau )\rightarrow 
\sigma \tau $, $\sigma ,\tau \in \FF _N^+$.
Let $\cP_N$ be the algebra of polynomials
on $N$ non-commuting
indeterminates $X_1$,$\ldots $,$X_N$ 
with complex coefficients. For any 
$\sigma=i_1\cdots i_l\in \FF_N^+$, 
we define
$X_\sigma=X_{i_1}\cdots X_{i_l}$. Using this notation, 
each element $P\in \cP_N$ can be uniquely written 
as
\begin{equation}\label{pesi}
P=\sum _{\sigma \in \FF _{N}^+}c_{\sigma }X_{\sigma },
\end{equation} 
with only finitely many coefficients $c_{\sigma }$
different from zero.
The length of the highest $\sigma $ such that $c_{\sigma }\ne 0$
is the {\it degree } of $P$.
We also have
$$P=\sum _{k\geq 0}P_k=
\sum _{k\geq 0}\sum _{|\sigma |=k}c_{\sigma }X_{\sigma },
$$
where each $P_k$ belongs to the vector space $\cL ^N_k$ of homogeneous
polynomials of degree $k\geq 0$ in $N$ variables $X_1$, $\ldots $,
$X_N$. The dimension of $\cL ^N_k$ is $N^k$. 
An involution $+$ 
on $\cP _N$ can be introduced as follows: $X^+_k=X_{k}$, $k=1,\ldots ,N$;
on monomials, 
$(X_{\sigma })^+=X_{I(\sigma )}$; in general,
 if $P$ has the representation as in \eqref{pesi}
then 
\begin{equation}P^+=\sum _{\sigma \in \FF _{N}^+}\overline{c}_{\sigma }
X_{\sigma }^+,
\end{equation}
and $\cP_N$ is a unital, associative $*$-algebra over $\CC $.
Let $\phi $ be a strictly positive functional on $\cP_N$,
that is, $\phi $ is a linear unital map on $\cP _N$ and
$\phi (P^+P)>0$
for every $P\in \cP _N-\{0\}$. The Gelfand-Naimark-Segal construction 
applied to
$\phi $ gives a Hilbert space $\cH _{\phi }$ such that 
$\{X_{\sigma }\}_{\sigma \in \FF ^+_N}$ is a  
linearly independent family in $\cH _{\phi }$. The 
Gram-Schmidt procedure gives a family $\{\varphi _{\alpha  }\}
_{\alpha \in \FF ^+_N}$ of polynomials 
such that 
\begin{equation}\label{bond1}
\varphi _{\alpha  }=
\sum _{\beta \preceq \alpha }a_{\alpha ,\beta }X_{\beta },
\quad a_{\alpha ,\alpha }>0;
\end{equation}
\begin{equation}\label{bond2}
\langle \varphi _{\alpha }, \varphi _{\beta }\rangle _{\phi }=
\delta _{\alpha ,\beta },
\quad \alpha ,\beta \in \FF ^+_N,
\end{equation}
where for $P_1, P_2\in \cP _N$,
$$\langle P_1,P_2\rangle _{\phi }=
\phi (P^+_2P_1).
$$
The elements $\varphi _{\alpha  }$, $\alpha \in \FF ^+_N$, will be called
the orthonormal polynomials associated with $\phi $.
We notice that the use of the Gram-Schmidt process depends
on the order chosen on 
$\FF _N^+$. A different order would give a different family of 
orthogonal polynomials. Due to the natural grading on 
$\FF _N^+$ it is possible to develop a base free approach to 
orthogonal polynomials. In the case of orthogonal polynomials
on several commuting variables this is presented in \cite{DX}. 
However, in this paper we consider only the lexicographic order on  
$\FF _N^+$.

The moments of $\phi $ are 
$$s_{\sigma }=\phi (X_{\sigma }), \quad \sigma \in \FF ^+_N,
$$ 
and we define the moment kernel of $\phi $ by the formula 
$K_{\phi }(\alpha ,\beta )=
s_{I(\alpha )\beta }$, $\alpha ,\beta \in \FF ^+_N$.
We notice that $K_{\phi }$ is a positive definite kernel
on $\FF ^+_N$ and for 
$\alpha ,\sigma ,\tau \in \FF _N^+$,
\begin{equation}\label{hankel}
K_{\phi }(\alpha \sigma ,\tau)=K_{\phi }(\sigma  ,I(\alpha )\tau ).
\end{equation}
This property can be viewed as a Hankel type condition.
Conversely, it is easily seen that if a positive definite kernel
$K$ satisfies \eqref{hankel} then there exists a positive
functional $\phi $ on $\cP _N$ such that $K=K_{\phi }$.

\bigskip
\noindent
{\em 2.1. Three term relations.}
Let $\phi $ be a unital strictly positive functional on 
$\cP _N$ and let $\{\varphi _{\alpha }\}_{\alpha \in \FF ^+_N}$
be the orthonormal polynomials associated with $\phi $.
As in the commutative case (see \cite{DX}),
it is very convenient to use a matrix notation, 
$\Phi _n=\left[\varphi _{\sigma }\right]_{|\sigma |=n}$
for $n\geq 0$ and $\Phi _{-1}=0$.
With this notation, the analogy with the classical case $N=1$ will be much 
more transparent. It turns out that the family $\{\Phi _n\}_{n\geq 0}$
satisfies a three-term recursive formula,
\begin{equation}\label{3termeni}
X_k\Phi _n=
\Phi _{n+1}A_{n+1,k}+
\Phi _nB_{n,k}+
\Phi _{n-1}A^*_{n,k}, 
\end{equation}
for $k=1,\ldots ,N$ and $n\geq 0$ (see \cite{Co} and \cite{BCJ}).
Each matrix $B_{n,k}$, $n\geq 0$, $k=1,\ldots ,N$, 
is a selfadjoint $N^n\times N^n$ matrix, while
each $A_{n,k}$, $n>0$, $k=1,\ldots ,N$, 
is an $N^{n}\times N^{n-1}$ matrix such that
$$A_n=\left[\begin{array}{ccc}
A_{n,1} & \ldots & A_{n,N}
\end{array}
\right]
$$
is an upper triangular invertible matrix for every $n\geq 0$,
with strictly positive elements on the diagonal.
The fact that $A_n$ is upper triangular comes from the lexicographic order 
that we use on 
$\FF ^+_N$. The invertibility of $B_n$ is a consequence of the fact
that $\phi $ is strictly positive and appears to be a basic translation of 
this information. The diagonal of $A_n$ is strictly 
positive since we chose $a_{\alpha ,\alpha }>0$. A 
family $\cA =\{A_{n,k}, B_{m,k}\mid n>0, m\geq 0, k=1,\ldots ,N\}$
of matrices satisfying all these properties will be called admissible.
It turns out that there are no other restrictions on the 
matrices $A_{n,k}$, $B_{n,k}$ as shown by the following Favard type 
result mentioned 
in \cite{Co}. A similar result for the monic orthogonal polynomials,
$p_{\sigma }=\displaystyle\frac{1}{a_{\sigma ,\sigma }}\varphi _{\sigma }$
was recently mentioned in \cite{An}.

\begin{theorem}\label{T4}
Let $\varphi _{\sigma }=\sum _{\tau \preceq \sigma }
a_{\sigma ,\tau}X_{\tau }$, $\sigma \in \FF _N^+$,
be elements in $\cP _N$ such that  $\varphi _{\emptyset }=1$ and 
$a_{\sigma ,\sigma }>0$.
Assume that there exists an admissible family 
$\cA $
of matrices such that 
for $k=1,\ldots ,N$ and $n\geq 0$, 
$$
X_k\Phi _n=
\Phi _{n+1}A_{n+1,k}+
\Phi _nB_{n,k}+
\Phi _{n-1}A^*_{n,k},
$$
where $\Phi _n=\left[\varphi _{\sigma }\right]_{|\sigma |=n}$
for $n\geq 0$ and $\Phi _{-1}=0$.
Then there exists a unique strictly positive functional $\phi $
on $\cP _N$ such that 
$\{\varphi _{\sigma }\}_{\sigma \in \FF _N^+}$
is the family of orthonormal  polynomials associated to $\phi $.
\end{theorem}

There is also a family of Jacobi matrices associated 
with the three-term relation
in the following way.  
For $P\in \cP _N$ define 
$$
\Psi _{\phi }(P)\varphi _{\sigma }=P\varphi _{\sigma }.
$$
Since the moment kernel has the Hankel type structure
mentioned above in \eqref{hankel}, it follows that each $\Psi _{\phi }(P)$ 
is a symmetric operator
on the Hilbert space $\cH _{\phi }$ with dense domain $\cP _N$. 
Moreover, 
for $P,Q \in \cP _N$, 
$$\Psi _{\phi }(PQ)=\Psi _{\phi }(P)\Psi _{\phi }(Q),$$
and $\Psi _{\phi }(P)\cD \subset \cD$,  
hence $\Psi _{\phi }$ is an unbounded representation of
$\cP _N$.
Also,  $\phi (P)=
\langle \Psi _{\phi }(P)1,1\rangle _{\phi }$ for $P\in \cP _N$.
We distinguish the operators $\Psi _k=\Psi _{\phi }(X_k)$,
$k=1,\ldots ,N$, since
$\Psi _{\phi }(\sum _{\sigma \in \FF ^+_N}c_{\sigma }X_{\sigma })
=\sum _{\sigma \in \FF ^+_N}c_{\sigma }\Psi _{\sigma }$. 
Let $\{e_1,\ldots ,e_N\}$ be the standard basis
of $\CC ^N$ and define the unitary operator $W$ from $l^2(\FF ^+_N)$
onto $\cH _{\phi }$ such that 
$W(e_{\sigma })=\varphi _{\sigma }$, $\sigma \in \FF ^+_N$.
We see that $W^{-1}\cD $ is the linear space $\cD _0$
generated by $e_{\sigma }$, $\sigma \in \FF ^+_N$, 
so that we can define
$$J_k=W^{-1}\Psi _{k}W,\quad k=1,\ldots ,N,$$ 
on $\cD _0$. Each $J_k$ is a symmetric operator on $\cD _0$
and by \eqref{3termeni}, the matrix of (the closure of) $J_k$ with respect
to the orthonormal basis 
$\{e_{\sigma }\}_{\sigma \in \FF ^+_N}$ is
$$J_k=\left[
\begin{array}{cccc}
B_{0,k} & A^*_{1,k} & 0 & \ldots \\
 & & & \\
A_{1,k} & B_{1,k} & A^*_{2,k} & \\
 & & & \\ 
0 & A_{2,k} & B_{2,k} & \ddots \\
 & & & \\ 
\vdots & & \ddots & \ddots 
\end{array}
\right].$$
We call $(J_1,\ldots ,J_N)$ a Jacobi $N$-family on $\cD _0$. 
It turns out that the usual admissibility conditions 
on $A_{n,k}$ and $B_{n,k}$ insure a
joint model of a Jacobi family in the following sense.
\begin{theorem}\label{jacobi}
Let $(J_1,\ldots ,J_N)$ be a Jacobi $N$-family 
and assume that the corresponding $\cA $
is an admissible family of matrices.
Then there exists a unique strictly positive
functional $\phi $ on $\cP _N$ with associated orthonormal
polynomials $\{\varphi _{\sigma }\}_{\sigma \in \FF ^+_N}$ such that the map
$W(e_{\sigma })=\varphi _{\sigma }$, $\sigma \in \FF ^+_N,$
extends to a unitary operator from $l^2(\FF ^+_N)$ onto $\cH _{\phi }$ and 
$J_k=W^{-1}\Psi _{k}W$ for $k=1,\ldots ,N$.
\end{theorem}
For details about the proof of this result see \cite{Co}.

\bigskip
\noindent
{\em 2.2. Jacobi $N$-families and combinatorics of lattice paths.}
The matrices $A_{n,k}$ and $B_{n,k}$ contain the whole information
about the orthonormal polynomials (or the moment kernel $K_{\phi }$).
Ususally they are called the Jacobi coefficients of $K_{\phi }$ 
and can be calculated from the moments. For instance, 
$$A_n=\left[a_{\alpha ,\beta }\right]^{-1}_{|\alpha |=|\beta |=n}
\left[a_{\alpha ,\beta }\right]^{\oplus N}_{|\alpha |=|\beta |=n-1},$$
where $a_{\alpha ,\beta }$ are the coefficients of the orthogonal
polynomials and 
for a matrix $A$ 
we use the notation 
$$A^{\oplus l}=
\underbrace{A\oplus \ldots \oplus A}_{\mbox{$l$ times}}.
$$
In their turn, the coefficients  $a_{\alpha ,\beta }$, $\beta \preceq 
\alpha $, can be calculated from the formula
$$a_{\alpha ,\beta }=\displaystyle\frac{1}
{\sqrt{D_{\alpha -1}D_{\alpha }}}
\det \left[K(\alpha ',\beta ')\right]_{\alpha '\prec \alpha ,
\beta '\preceq \alpha ,\beta '\ne \beta },$$
where 
$$D_{\alpha }=\det \left[K(\alpha ',\beta ')\right]_{\alpha ',\beta 
\preceq \alpha }.$$
The formula for $B_{n,k}$ is somewhat more involved and we do not 
record it here. Instead we consider a different kind of relation between 
moments and Jacobi coefficients which appears to be more explicit.
The case $N=1$ 
is classical, see \cite{Fl}, \cite{La}.
The Jacobi $N$-family $\left(J_1, \ldots ,J_N\right)$ 
is a convenient tool to deal with this matter.
Thus, for any $\sigma \in \FF ^+_N$ we have that 
$$s_{\sigma }=\phi (X_{\sigma })=
\langle \Psi _{\phi }(X_{\sigma })1,1 \rangle _{\phi }=
\langle \Psi _{\sigma }1,1\rangle _{\phi },
$$
and by Theorem ~\ref{jacobi}, 
$$\langle \Psi _{\sigma }1,1\rangle _{\phi }=
\langle J_{\sigma }e_0,e_0\rangle ,
$$
therefore we have 
\begin{equation}\label{key}
s_{\sigma }=\langle J_{\sigma }e_0,e_0\rangle ,\quad \sigma \in \FF ^+_N.
\end{equation}

Now we introduce some special paths on  
$\NN _0\times \{1,\ldots ,N\}\times \NN _0$. The allowed steps 
are the following: level steps $l^k_{n,m}$ from a point 
$(n,k,m)$ to $(n+1,k,m)$, level steps $lp^n_{m,k}$
from a point $(n,k,m)$ to $(n,p,m)$, for some $p\in \{1,\ldots ,N\}-\{k\}$,
rise steps $r^k_{n,m}$ from   
a point 
$(n,k,m)$ to $(n+1,k,m+1)$, 
and fall steps $f^k_{n,m}$ from a point 
$(n,k,m)$ to $(n+1,k,m-1)$ (see Figure ~1 for an example).

\begin{figure}[h]
\setlength{\unitlength}{3000sp}%
\begingroup\makeatletter\ifx\SetFigFont\undefined%
\gdef\SetFigFont#1#2#3#4#5{%
  \reset@font\fontsize{#1}{#2pt}%
  \fontfamily{#3}\fontseries{#4}\fontshape{#5}%
  \selectfont}%
\fi\endgroup%
\begin{picture}(3399,2499)(214,-1948)
\thinlines
{\put(601,-1936){\vector( 0, 1){2475}}
}%
{\put(226,-1561){\vector( 1, 0){2775}}
}%
{\multiput(976,-1261)(9.01163,0.00000){259}{\makebox(1.6667,11.6667){\SetFigFont{5}{6}{\rmdefault}{\mddefault}{\updefault}.}}
}%
{\multiput(1426,-961)(8.98760,0.00000){243}{\makebox(1.6667,11.6667){\SetFigFont{5}{6}{\rmdefault}{\mddefault}{\updefault}.}}
}%
{\put(601,-1561){\vector( 2, 1){1500}}
}%
\thicklines
{\put(1201,-1261){\line( 1, 0){300}}
\put(1501,-1261){\line( 2, 3){300}}
\put(1801,-811){\line( 2, 1){600}}
\put(2401,-511){\line( 2,-3){300}}
\put(2701,-961){\line( 1, 0){300}}
}%
\end{picture}%
\caption{\mbox{An example of a path for $N=2$}}
\end{figure}

We introduce a weight on steps by the formula
$$
w(\mbox{step})=\left\{
\begin{array}{cl}
I & \mbox{if step $=lp^n_{m,k}$}\\
B_{m,k} & \mbox{if step $=l^k_{n,m}$}\\
A_{m+1,k} & \mbox{if step $=r^k_{n,m}$}\\
A^*_{m,k} & \mbox{if step $=f^k_{n,m}$},
\end{array}
\right.
$$
where $I$ denotes the identity matrix of appropriate size.
If ${\bf p}$ is made of $l$ steps, step $1$, $\ldots $, step $l$, 
then we define the weight of ${\bf p}$ by the formula
$$w({\bf p})=w(\mbox{step $l$})\ldots w(\mbox{step $1$}).$$
Any word 
$\sigma \in \FF ^+_N-\{\emptyset \}$ has a unique representation
$\sigma =i_1^{k_1}\ldots i_p^{k_p}$ with 
$1_1$, $\ldots $, $i_p\in \{1, \ldots ,N\}$, $k_1$, $\ldots $, $k_p>0$, and 
$i_l\ne i_{l+1}$ for $l=1$, $\ldots $, $p-1$. We consider 
the set $\cM _{\sigma }$ of all paths that start at $(0,i_p,0)$ and 
end at $(|\sigma |,i_1,0)$, with the property that 
the first $k_p$ steps belong to $\NN _0\times \{i_p\}\times \NN _0$, 
the next $k_{p-1}$ steps belong to $\NN _0\times \{i_{p-1}\}\times \NN _0$, 
and so on, until the last $k_1$ steps which belong to 
$\NN _0\times \{i_1\}\times \NN _0$.
These sets are related to the set of Motzkin paths.
The Motzkin paths of length $n$ are the 
paths in $\NN ^2_0$ 
made of level, fall, and rise steps, starting at  
$(0,0)$  and 
ending at $(n,0)$. Their set is denoted by $\cM _n$ and the number of 
elements of $\cM _n$ is given by the Motzkin number
$$M_n=\frac{1}{n}\sum _k \left(\begin{array}{c} n \\ k \end{array}\right)
\left(\begin{array}{c} n-k \\ k-1 \end{array}\right).
$$
It is easily seen that for any $\sigma \in \FF ^+_N-\{\emptyset \}$
there is a bijection between 
$\cM_{\sigma }$ and $\cM _{|\sigma |}$.  
We can now describe a combinatorial structure of the 
moments.
\begin{theorem}\label{comb}
Let  $\phi $
be a  strictly positive
functional on $\cP _N$
and let $\cA$
be the admissible family of matrices associated with 
$\phi $ by \eqref{3termeni}.
Then the moments of $\phi $ can be calculated by the formula
\begin{equation}\label{momjac}
s_{\sigma }=\sum _{{\bf p}\in \cM _{\sigma }}w({\bf p}),
\quad \quad \sigma \in \FF ^+_N-\{\emptyset \}.
\end{equation}
\end{theorem}
\begin{proof}
We consider the following points in 
$\NN _{0}\times \{1,\ldots ,N\}\times \NN _{0}$:
$$P_{n,j}=(0,j,n), \quad n\geq 0, j=1,\ldots ,N;
$$
$$Q_{n,k,m}=(n,k,m), \quad m,n\geq 0, k=1,\ldots ,N.
$$
For $\sigma \in \FF ^+_N-\{\emptyset \}$, 
$\sigma =i_1^{k_1}\ldots i_p^{k_p}$, we claim that 
$$J_{\sigma }=\left[
\begin{array}{ccc}
J^{\sigma }_{0,0} & J^{\sigma }_{0,1} & \ldots \\
J^{\sigma }_{1,0} & J^{\sigma }_{1,1} &    \\
\vdots & &\ddots 
\end{array}
\right],$$
where the entry $J^{\sigma }_{k,j}$
gives the sum of (the weights of) the paths in 
$\NN _{0}\times \{1,\ldots ,N\}\times \NN _{0}$ from
$P_{j,i_p}$ to $Q_{|\sigma |,i_1,k}$.
The claim is clearly true for $|\sigma |=1$ and then suppose it true for any 
word of length $\leq n$. 
Then consider a word $\sigma $ of length $n+1$. Several cases can occur.

{\it Case 1.} $k=0$. Let $\sigma =i_1^{k_1}\ldots i_p^{k_p}=i_1\tau $.
First, assume $k_1=1$. Due to the fact that the level steps of type 
$lp$ have weight $I$ and by the induction hypothesis, we deduce that 
the sum of the paths from 
$P_{j,i_p}$ to $Q_{|\sigma |,i_1,0}$
is 
$$B_{0,i_1}J^{\tau }_{0,j}+A^*_{1,i_1}J^{\tau }_{1,j},$$ 
which is precisely the $(0,j)$ entry of the product
$$
J_{i_1}J_{\tau }=J_{\sigma }.
$$
The case $k_1>1$ is similar, just by the induction hypothesis we deduce
that the sum of paths from 
$P_{j,i_p}$ to $Q_{|\sigma |,i_1,0}$
is again 
$$B_{0,i_1}J^{\tau }_{0,j}+A^*_{1,i_1}J^{\tau }_{1,j},$$ 
which is precisely the $(0,j)$ entry of the product
$$
J_{i_1}J_{\tau }=J_{\sigma }.
$$

{\it Case 2.} $j=0$ is similar.

{\it Case 3.} $k,j\geq 1$, then the induction hypothesis
implies that the sum of paths from 
$P_{j,i_p}$ to $Q_{|\sigma |,i_1,k}$
is
$$A_{k,i_1}J^{\tau }_{k-1,j}+B_{k,i_1}J^{\tau }_{k,j}+
A^*_{k+1,i_1}J^{\tau }_{k+1,j},
$$
which is precisely the $(k,j)$ entry of $J_{\sigma }$.
\end{proof}
When all $B_{n,k}$ are zero, the level steps of type $l^k_{n,m}$
dissapear, and our discussion is somewhat related to parts of \cite{Ni}.

Formula \eqref{momjac}
can be also used to calculate the Jacobi coefficients from the moments 
in a relatively
simple way. Let 
$\sigma \in \FF ^+_N-\{\emptyset \}$. It is convenient 
to introduce the notation $i(n)$ in order to denote
the $n$th letter of the word $\sigma $ (from left to right). 
If $|\sigma |=2n$, then there exists a unique path ${\bf p}_{\sigma }$
with corresponding weight
$$A^*_{1,i(i)}\ldots A^*_{n,i(n)}A_{n,i(n)}\ldots A_{1,i(i)}.$$
If $|\sigma |=2n+1$, then there exists a unique path, 
still denoted ${\bf p}_{\sigma }$,
with corresponding weight
$$A^*_{1,i(i)}\ldots A^*_{n,i(n)}B_{n+1,i(n+1)}A_{n,i(n)}\ldots A_{1,i(i)}.$$ 
In any case, let $\cM ^*_{\sigma }=\cM _{\sigma }-\{{\bf p}_{\sigma }\}$.
Also, we introduce the notation: $\tilde A_1=A_1$, and for $n\geq 2$,
$$\tilde A_n=A_nA^{\oplus N}_{n-1}\ldots A^{\oplus N^{n-1}}_{1}.$$
\begin{corollary}\label{jacmom}
The following formulae hold: for $n\geq 1$,
$$A^*_nA_n=\left(\tilde A^*_{n-1}\right)^{-1}
\left(\left[K_{\phi }(\sigma ,\tau )\right]
_{|\sigma |=|\tau |=n}-
\left[\sum _{{\bf p}\in \cM ^*_{I(\sigma )\tau }}w({\bf p})\right]_
{|\sigma |=|\tau |=n}
\right)\tilde A_{n-1}^{-1};$$
$$B_{0,k}=s_k, \quad k=1,\ldots ,N,$$
and 
for $n\geq 1$, $k=1,\ldots ,N,$
$$B_{n,k}==\left(\tilde A^*_{n}\right)^{-1}
\left(\left[K_{\phi }(k\sigma ,\tau )\right]
_{|\tau |=|\sigma |+2=n+1}-
\left[\sum _{{\bf p}\in \cM ^*_{I(\sigma )\tau }}w({\bf p})\right]_
{|\tau |=|\sigma |+2=n+1}
\right)\tilde A_{n}^{-1}.$$
\end{corollary}
Due to the fact that $A_n$ is an upper triangular matrix with strictly 
positive elements on the diagonal, the first relation of the previous 
result uniquely determine $A_n$ by Cholesky factorization.

\section{Free products}
The set $\cP _N$ can be viewed as the free product of $N$ 
copies of $\cP _1$:
$$\cP _N=\underbrace{\cP _1\star \ldots \star \cP _1}_{\mbox{$N$ times}}=
\CC\oplus 
\left(\oplus _{n\geq 1}
\oplus _{i_1\ne i_2,\ldots ,i_{n-1}\ne i_n}\cP ^0_{i_1}
\otimes \ldots \otimes \cP ^0_{i_n}\right),$$
where $\cP ^0_{i}$ is the set of polynomials in the variable $X_i$, 
$i=1,\ldots ,N$, without constant term.
This remark suggests that the simplest examples of families of 
orthogonal polynomials can be obtained by using free products. Some examples 
already appeared in \cite{An}. Here we describe a genreal construction. 
This allows to introduce multivariable analogues of all 
classical orthogonal polynomials.

The simplest attempt to construct families of orthogonal polynomials
on $\cP _N$ would be to consider orthogonal polynomials associated
with free products of positive functionals. Let $\phi _1 $ 
and 
$\phi _2 $ be two strictly positive functional on 
$\cP _{N_1}$, respectively $\cP _{N_2}$. Their free
product $\phi =\phi _1\star \phi _2$ on $\cP _{N_1}\star \cP _{N_2}$
is defined by $\phi (1)=1$ and 
$\phi (P_{i_1}\ldots P_{i_n})=\phi _{i_1}(P_{i_1})\ldots 
\phi _{i_n}(P_{i_n})$ for $n\geq 1$, 
$i_1\ne i_2$, $\ldots $, $i_{n-1}\ne i_n$, $P_{i_k}\in \cP ^0_{N_{i_k}}$,
and $i_k\in \{1,2\}$ for $k=1,\ldots ,n$.  
By results in \cite{Boc}, \cite{Boz}, $\phi $ is a positive functional. 
However, as it turns out, $\phi $ is not strictly positive. Thus, consider
the case of two strictly positive functionals  
$\phi _1 $ 
and 
$\phi _2 $ on $\cP _1$. Let $K$ be the restriction of the 
moment kernel $K_{\phi _1\star \phi _2}$ to the set 
$\{1, 12\}$. Its matrix is then 
$$
\begin{array}{rcl}
K& =&\left[
\begin{array}{cc}
\phi _1(X^2_1) & \phi _1(X^2_1)\phi _2(X_2) \\
\phi _2(X_2)\phi _1(X^2_1) & \phi _2(X_2)\phi _1(X^2_1)\phi _2(X_2)
\end{array}
\right] \\
 & & \\
& =& 
\left[\begin{array}{cc}
1 & 0 \\
0 & \phi _2(X_2)
\end{array}
\right]
\left[\begin{array}{cc}
\phi _1(X^2_1) & \phi _1(X^2_1) \\
\phi _1(X^2_1) & \phi _1(X^2_1)
\end{array}
\right]
\left[\begin{array}{cc}
1 & 0 \\
0 & \phi _2(X_2)
\end{array}
\right].
\end{array}
$$
The matrix $\left[\begin{array}{cc}
1 & 1 \\
1 & 1
\end{array}
\right]$ has rank one, so $K$ is never invertible.

We can consider another simple way related to free products 
in order to build orthogonal 
polynomials in several noncommutative variables. Thus, 
let $\{\varphi _{n,k}\}$, $n\geq 0$, $k\in \{1,\ldots ,N\}$, be $N$
families of orthonormal polynomials on the real line, determined
by the recursion formulae:
\begin{equation}\label{onedim}
x\varphi _{n,k}(x)=a_{n+1,k}\varphi _{n+1,k}(x)+
b_{n,k}\varphi _{n,k}(x)+a_{n,k}\varphi _{n-1,k}(x).
\end{equation}
We introduce polynomials in $N$ noncommutative variables 
as follows. Any word $\sigma \in \FF ^+_N-\{\emptyset \}$ can be uniquely
represented in the form $\sigma =i^{k_1}_1\ldots 
i^{k_p}_p,$ 
$i_1$, $\ldots $, $i_p\in \{1, \ldots ,N\}$, $k_1$, $\ldots $, $k_p>0$, and 
$i_l\ne i_{l+1}$ for $l=1$, $\ldots $, $p-1$.
Then define
\begin{equation}\label{freeprod}
\varphi _{\sigma }(X_1,\ldots ,X_N)=\varphi _{k_1,i_1}(X_{i_1})\ldots 
\varphi _{k_p,i_p}(X_{i_p}).
\end{equation}
\begin{theorem}\label{main}
There exists an admissible family $\cA $ of matrices such that 
for $k=1,\ldots ,N$ and $n\geq 0$,
\begin{equation}\label{33termeni}
X_k\Phi _n=
\Phi _{n+1}A_{n+1,k}+
\Phi _nB_{n,k}+
\Phi _{n-1}A^*_{n,k}, 
\end{equation}
where  $\Phi _n=\left[\varphi _{\sigma }\right]_{|\sigma |=n}$
for $n\geq 0$, $\Phi _{-1}=0$, and 
$\varphi _{\sigma }$ are given by \eqref{freeprod}.
\end{theorem}
\begin{proof}
First we prove the result for $N=2$.
From \eqref{onedim} we deduce
$$X_1=\varphi _1a_{1,1}+b_{0,1}=
\left[
\begin{array}{cc}
\varphi _1 & \varphi _2
\end{array}
\right]
\left[
\begin{array}{c}
a_{1,1}\\
0 
\end{array}
\right]
+b_{0,1}
$$
so that 
$$X_1\Phi _0=\Phi _1A_{1,1}+\Phi _0B_{0,1}
$$
with 
$$A_{1,1}=\left[
\begin{array}{c}
a_{1,1} \\
0
\end{array}
\right]\quad  \mbox{and}\quad  
B_{0,1}=\left[
\begin{array}{c}
b_{0,1}
\end{array}
\right].
$$
Similarly, 
$$X_2\Phi _0=\Phi _1A_{1,2}+\Phi _0B_{0,2}
$$
with 
$$
A_{1,2}=\left[
\begin{array}{c}
a_{1,2} \\
0
\end{array}
\right]\quad  \mbox{and}\quad  
B_{0,2}=\left[
\begin{array}{c}
b_{0,2}
\end{array}
\right].
$$
We see that 
$$A_1=
\left[
\begin{array}{cc}
A_{1,1} & A_{1,2}
\end{array}
\right]
=
\left[
\begin{array}{cc}
a_{1,1} & 0 \\
0 & a_{1,2}
\end{array}
\right]
$$
is upper triangular (actually diagonal) and has the 
elements on the diagonal $>0$ (since all $a_{n,k}>0$).
The case $n\geq 1$ can be delt with in a similar manner. 
The first $2^{n-1}$ words of length
$n$ start with letter $1$ and have the structure
$$1^{n-k}\tau $$
where $\tau $ is a word of lenght $k$ starting with letter $2$
(unless it is $\emptyset $). 
For $k=0$ there is exactly one such word,
$1^n$, while for $0<k\leq n-1$, there 
are $2^{k-1}$ such kind of words. Using \eqref{onedim}
we have that 
$$
\begin{array}{rcl}
X_1\varphi _{1^{n-k}\tau }&=&
X_1\varphi _{1^{n-k}}\varphi _{\tau }=X_1\varphi _{n-k,1}\varphi _{\tau }\\
& & \\ 
&=&\varphi _{n-k+1,1}\varphi _{\tau }a_{n-k+1,1}+
\varphi _{n-k,1}\varphi _{\tau }b_{n-k,1}+
\varphi _{n-k-1,1}\varphi _{\tau }a_{n-k,1} \\
& & \\ 
&=&\varphi _{1^{n-k+1}\tau }a_{n-k+1,1}+
\varphi _{1^{n-k}\tau }b_{n-k,1}+
\varphi _{1^{n-k-1}\tau }a_{n-k,1}.
\end{array}
$$
The last $2^{n-1}$ words of lenght $n$ start with letter $2$
and therefore we have for such a word $\sigma $ that
$$
\begin{array}{rcl}
X_1\varphi _{\sigma }&=&X_1\varphi _{0,1}\varphi _{\sigma } \\
 & & \\
 &=&\varphi _{1,1}\varphi _{\sigma }a_{1,1}+\varphi _{\sigma }b_{0,1} \\
 & & \\
 &=&\varphi _{1\sigma }a_{1,1}+\varphi _{\sigma }b_{0,1}.
\end{array}
$$
Putting together the above two formulae we deduce that
$$
X_1\Phi _n=
\Phi _{n+1}A_{n+1,1}+
\Phi _nB_{n,1}+
\Phi _{n-1}A^*_{n,1}, 
$$
where 
$$A_{n,1}=
\left[
\begin{array}{ccclc}
a_{n,1} & & & & \\
 & a_{n-1,1} & & & \\
 & & a^{\oplus 2}_{n-2,1} & & \\
 & & & \ddots & \\
 & & & & a^{\oplus 2^{n-2}}_{1,1} \\
 & & & & \\
 & & 0_{2^{n-1}\times 2^{n-1}} & & 
\end{array}
\right]
$$
and 
$$B_{n,1}=
\left[
\begin{array}{ccclc}
b_{n,1} & & & & \\
 & b_{n-1,1} & & & \\
 & & b^{\oplus 2}_{n-2,1} & & \\
 & & & \ddots & \\
 & & & & b^{\oplus 2^{n-2}}_{1,1} \\
\end{array}
\right];
$$
the unspecified entries are all zero.
Similarly, we deduce 
$$
X_2\Phi _n=
\Phi _{n+1}A_{n+1,2}+
\Phi _nB_{n,1}+
\Phi _{n-1}A^*_{n,2}, 
$$
for some suitable matrices $A_{n,2}$ and $B_{n,2}$. Actually, the same 
proof works for $N>2$ and we record here the form of the 
matrices $A_{n,k}$, $B_{n,k}$ for an arbitrary $N$. Let 
$$\cW ^n_k=\{\sigma \in \FF ^+_N\mid \mbox{$|\sigma |=n$
and $\sigma =k\tau $ for some $\tau $}\}$$
and denote by $\pi ^n_k$ the bijection from the set 
$\{1,2,\ldots ,N^{n-1}\}$ onto $\cW ^n_k$ defined simply by 
$\pi ^n_k(l)=$the $l$th word in $\cW ^n_k$, with respect to the 
lexicographic order. 
Also, if $\sigma \in \cW ^n_k$ then it has a unique representation 
$\sigma =k^p\tau $ with $\tau $ a word that does not start with the 
letter $k$. Define $n_k(\sigma )=p$. 
Now $B_{n,k}$ is an $N^n\times N^n$
matrix such that, for $l,m\in \{1, 2,\ldots ,N^n\}$,
\begin{equation}\label{b}
\left(B_{n,k}\right)_{m,l}=
\left\{\begin{array}{cl}
b_{n_k(\pi ^{n+1}_k(m))-1,k} & l=m \\
 &\\
0 & l\ne m
\end{array}
\right.
\end{equation}
and $A_{n,k}$ is an $N^n\times N^{n-1}$ matrix such that 
for $l\in \{1,2,\ldots ,N^{n-1}\}$ and $m \in \{1, 2,\ldots ,N^n\}$,
\begin{equation}\label{a}
\left(A_{n,k}\right)_{m,l}=
\left\{\begin{array}{cl}
a_{n_k(\pi ^{n}_k(m))-1,k} & l=m \\
 & \\
0 & l\ne m
\end{array}.
\right.
\end{equation}
Therefore $A_n=
\left[
\begin{array}{ccc}
A_{n,1}
& \ldots & A_{n,N}
\end{array}
\right]
$ 
is a diagonal matrix with strictly positive diagonal elements, so that 
$\cA =\{A_{n,k}, B_{m,k}\mid n>0, m\geq 0, k=1,\ldots ,N\}$ is an 
admissible family of matrices.
\end{proof}

\end{document}